\documentstyle[epsf]{article}
\newtheorem{definition}{Definition}
\newtheorem{theorem}{Theorem}

\newtheorem{lemma}{Lemma}
\newtheorem{proposition}{Proposition}
\newtheorem{corollary}{Corollary}

\newfont{\germ}{eufm10}
\def\goth#1{\mbox{\germ #1}}

\def\A{{\cal A}}
\def\bi{\goth{b}}

\def\ch{\mbox{\sl ch}}
\def\en{\goth{n}}
\def\ep{\epsilon}
\def\et#1{\tilde{e}_{#1}}
\def\ft#1{\tilde{f}_{#1}}
\def\ha{\goth{h}}

\def\La{\Lambda}
\def\la{\lambda}
\def\P{{\cal P}}
\def\Proof{\noindent{\sl Proof.}\quad}
\def\Pset{\P_L(\La)}
\def\qed{~\rule{1mm}{2.5mm}}
\def\slth{\widehat{sl}(2)}
\def\T{{\cal T}}
\def\Uq#1{U_q\bigl(#1\bigr)}
\def\vep{\varepsilon}
\def\wt{\mbox{\sl wt}~}
\def\Y{{\bf Y}}
\def\Yset{{\cal Y}(\La)}
\def\Z{\:\mbox{\sf Z} \hspace{-0.82em} \mbox{\sf Z}\,}
\def\Zn{\Z_{\ge0}}
\def\Zs{\mbox{\scriptsize \rm Z}  \! \! \mbox{\scriptsize \rm Z}}
\def\Zss{\mbox{\tiny \rm Z}  \! \! \mbox{\tiny \rm Z}}

\begin{document}
\title{ Demazure modules and vertex models:\\ 
        the $\slth$ case }

\author{Omar Foda$^1$\thanks{Supported in part by the Australian Research 
                             Council (ARC), and the Netherlands Organization 
                             for Scientific Research (NWO).}, 
        Kailash C. Misra$^2$\thanks{Supported in part by NSA/MSP Grant 
                                    No. MDA92-H-3076.}
        and Masato Okado$^3$\thanks{Supported in part by Grant-in-Aid for 
                                    Scientific Research on Priority Areas, 
                                    the Ministry of Education, Science and 
                                    Culture, Japan}}

\date{ \it $^1$\it Instituut voor Theoretische Fysica,\\
       \it Universiteit Utrecht,\\
       \it Utrecht 3508 TA, The Netherlands \thanks{ 
       Permanent address: Department of Mathematics, 
       The University of Melbourne, Parkville, Victoria 
       3052, Australia}\\
       \vskip 4mm
       \it $^2$Department of Mathematics,\\
       \it North Carolina State University,\\
       \it Raleigh, NC 27695-8205, USA\\ 
       \vskip 4mm
       \it $^3$Department of Mathematical Sciences,\\
       \it Faculty of Engineering Science,\\
       \it Osaka University, Toyonaka,\\
       \it Osaka 560, Japan }

\maketitle

\begin{abstract}
We characterize, in the case of $\slth$, the crystal base of the Demazure 
module $E_w(\La)$ in terms of extended Young diagrams or paths for any dominant 
integral weight $\La$ and Weyl group element $w$. Its character is evaluated 
via two expressions, 'bosonic' and 'fermionic'.
\end{abstract} 

\setcounter{section}{-1}

\section{Introduction}

The purpose of this work is to characterize the crystal bases of the 
Demazure modules of $\slth$ in terms of 'paths', or extended Young
diagrams, and obtain explicit expressions for their full characters. 
There are two ways to motivate this work: one can motivate it from 
the viewpoint of certain recent developments in mathematical physics, 
or from the viewpoint of representation theory. 

\subsection{1-point functions in exactly solvable models}

Certain physical quantities in exactly solvable two-dimensional 
lattice models, namely the so-called '1-point functions' can be 
evaluated using Baxter's corner transfer matrix method.\footnote{For 
an introduction to exactly solvable models, and to the corner transfer
matrix method, we refer to \cite{BaxterBook}.} This method reduces 
the computation of 1-point functions to a computation of a weighted 
sum. This sum is over combinatorial objects called 'paths': 
1-dimensional configurations, defined on the half-line. They will 
be discussed in detail in the sequel. The weight of a path is the 
evaluation, on that path, of an 'energy' functional defined on the 
set of all paths. The weighted sum is the generating function of 
the number of paths of a certain weight.

One way to evaluate this generating function is by solving a recursion
relation for the generating function that counts the number of paths 
defined on a finite segment, of length $L$, of the half-line. The recursion 
is with respect to $L$. The $L \rightarrow \infty$ limit of the 
result, which is the generating function of all paths, turns out to be 
the character of a highest weight module of an infinite-dimensional 
algebra: an affine, or Virasoro algebra.\footnote{For
an introduction to the algebaric approach to exactly solvable models, we
refer to \cite{JimboMiwaBook}.} The result obtained by solving a recursion 
relation typically has a form that, for physical reasons, is called 
'bosonic'.\footnote{For an explanation of the origin of this terminology, 
we refer to \cite{StonyBrookReview}.} 

\subsection{Rogers-Ramanujan-type $q$-series identities}

More recently, the same physical objects were evaluated in a completely 
different way, using Bethe Ansatz methods \cite{StonyBrookReview}. In 
that case, the results turn out to have a completely different form, 
that, for the same physical reasons as above, is refered to as 'fermionic'. 

Equating the bosonic and fermionic expressions for the same
objects, one obtains Rogers-Ramanujan-type identities between
$q$-series. Strictly speaking, such identities are conjectures,
since the methods used to obtain them are indirect, and involve
certain 'physical assumptions'.\footnote{What we have in mind
are the various assumptions involved in the derivation of the 
corner transfer matrix method, and the 'string hypothesis'.}
What one needs are direct proofs. 

\subsection{Schur-type polynomial identities}

One method to obtain direct proofs, which dates back to Schur, 
is to work at the level of the generating functions of all paths
that live on a segement of the half-line of length $L$. 
These are polynomials that depend on $L$. Once we can prove the 
boson-fermion identities for finite $L$, we take the limit 
$L \rightarrow \infty$, and obtain the original $q$-series identities.

The following question now arises: If the original infinite
$q$-series are characters of highest weight modules of
infinite-dimensional algebras, are the finite-$L$ polynomials the
characters of anything? Do they have a meaning in mathematics? Or
are they just convenient objects that appear in intermediate
steps?

The reason why this question arises is that rigorous proofs in
exactly solvable models are typically obstructed by the fact that
we have to deal with infinite series, and infinite-dimensional
quantities. The reason why we typically have to do that, is that
only in the limit $L \rightarrow \infty$ do these model exhibit
invariance under infinite-dimensional algebras, and this
 invariance is an essential ingredient in solvability. Thus we
have the following problem: we need to work in the $L \rightarrow
\infty$ limit in order to be able to obtain exact solutions, but
in that limit we cannot provide rigorous proofs.

This is the reason why it is very interesting to try to
understand as much of the structure of exactly solvable models as
possible at the level of finite $L$. The hope is that we can find
algebraic structures at finite $L$ that are sufficiently strong
to provide rigorous proofs and solutions. A typical object to
investigate are the $L$-restricted generating functions.

What we find in this work, is that these objects do have a
mathematical significance, or are directly related to objects
that do. We find that the $L$-restricted generating functions are
closely related to the characters of Demazure modules.

\subsection{Demazure modules}

{}From the viewpoint of representation theory, we can motivate our 
work in the evaluation of the characters of Demazure modules. The 
Demazure module is characterized by its highest weight $\La$ and
Weyl group element $w$, and is denoted by $E_w(\La)$. 
In \cite{Ku,Ma}, a formula for computing the Demazure character 
$\ch E_w(\La)=\sum_{\mu\in P} \dim \bigl(E_w(\La)\bigr)_\mu e^\mu$ 
is given as follows. For $\mu\in P, i=0,1$, define the operator
$D_i:\Z[P]\longrightarrow\Z[P]$ by 
\[
D_i(e^\mu)=\frac{e^{\mu+\rho}-e^{r_i(\mu+\rho)}}{1-e^{-\alpha_i}}e^{-\rho}.
\]
Let $w=r_{i_n}\cdots r_{i_2}r_{i_1}\in W$ be a reduced expression. 
For the terminologies of $\slth$ see subsection 1.1. Then, the
Demazure character formula states that
\[
\ch E_w(\La)=D_{i_n}\cdots D_{i_2}D_{i_1}(e^\La).
\]
As elegant as this character formula is, it is anything but combinatorial.
Hence using this formula it is very difficult to get any information about
specific weight spaces. In \cite{S}, Sanderson used Littelrmann's path 
model \cite{Li} and gave nice expressions for the `real' characters of 
$E_w(\La)$ for all $\La$ and the `principal' characters of $E_w(\La)$ for
$\La=s\La_0$. We point out that the real (resp. principal) character of
$E_w(\La)$ is the specialization of $e^{-\La}\ch E_w(\La)$ where
$e^{-\alpha_0}=q,e^{-\delta}=1$ (resp. $e^{-\alpha_0}=q,e^{-\delta}=q^2$).

\subsection{Crystal bases}

In 1990 Kashiwara brought a notion of crystal base in the representation 
theory of the quantum group \cite{K0}. This notion is quite powerful in 
the combinatorial aspect of representation theory. As an example, we find 
\cite{Nakashima} in which the Littlewood-Richardson type rules for the 
classical Lie algebras are given in a purely combinatorial way. Kashiwara
also obtained the crystal bases of the Demazure modules as subsets of those 
of the corresponding highest weight modules \cite{K2}. His algorithm is 
given in a recursive way. Consequently he got a new proof of the Demazure
character formula. 

Since we know that in the case of affine Lie algebra the crystal base of the 
highest weight module is described in terms of paths \cite{(KMN)s1,(KMN)s2},
it is natural to try to characterize the crystal base of the Demazure module
in this language. For the affine algebra of type A we have yet another notion
to describe the crystal base, that is, extended Young diagram \cite{MM,JMMO}.

\subsection{Main results}

What we have done in our work is the characterization of Demazure crystal
bases in terms of extended Young diagrams (Theorem \ref{Demcry}). One of 
the advantages to do so is that we can make the most of earlier efforts
toward the evaluations of so called 1-dimensional configuration sums 
\cite{DJKMO,DJKMO1d}. Recent results on fermionic expressions are also
available. We would like to note here that after we had obtained the
fermionic expression (Theorem \ref{th3}), we came to know Schilling 
also had the same result \cite{Schilling}. 
However, we have included our proof since it is different from that of 
Schilling. As a corollary
of the characterization theorem, we evaluate the characters of Demazure
modules (Theorem \ref{fermionic} and \ref{Demch}).

\subsection{Outline of paper}

In section 1, we define the basic combinatorial objects in this work: 
paths and extended Young diagrams. We then describe the crystal structure
of the integrable highest weight module of $\Uq{\slth}$. In section 2, 
we introduce the Demazure modules, and characterize their crystals via
paths or extended Young diagrams. In section 3, we review 1-dimensional
configuration sums and rewrite them into fermionic forms. Using them
we can express Demazure characters. We also discuss some specializations
and compare them with Sanderson's results. Section 4 contains a short 
discussion.

\section{Paths and extended Young diagrams} 

In this section we recall the basic combinatorial objects used in
this work: the paths, and extended Young diagrams. Next, we outline 
the realization of the crystal base of an integrable $\slth$-module 
in terms of these objects. We also list some properties of the 
crystal base which we will need in the sequel.

\subsection{Preliminaries}
Consider the affine Lie algebra $\slth$ \cite{Kac}. 
Let 
$\{\alpha_0,\alpha_1\}$, $\{h_0,h_1\}$ and $\{\La_0,\La_1\}$ 
be the simple roots, simple co-roots and fundamental weights 
respectively.
They satisfy
$\La_i(h_j)=\delta_{ij}$ and 
$\alpha_i(h_j)=2(-1)^{1+\delta_{ij}}$ 
for $i,j=0,1$. 
$\delta=\alpha_0+\alpha_1$, $c=h_0+h_1$ 
are the null root and canonical central element, respectively. 
The sets
$P=\Z\La_0\oplus\Z\La_1\oplus\Z\delta$, and 
$Q=\Z\alpha_0\oplus\Z\alpha_1$
are the weight and root lattices, respectively. 
Let $\rho\in\ha^*$ be such that $\rho(h_i)=1$ for $i=0,1$.

Let $\Uq{\slth}$ denote the quantized universal enveloping algebra
associated with $\slth$. For its precise definition and Hopf algebra
structure, we refer to \cite{JMMO}.
For $\La\in P$,
$\La(c)$ is the level of $\La$. The set $P^+=\{\La\in P\mid \La(h_i)
\ge0,i=0,1\}$ is the set of dominant weights. 
For $\La\in P^+$, let $V(\La)$ denote the unique (up to isomorphism) 
integrable highest weight $\Uq{\slth}$-module.
Since $V(\La+k\delta)\simeq V(\La)\otimes V(k\delta)$ and 
$\dim V(k\delta)=1$, it suffices to assume that $\La=s\La_0+t\La_1$ 
for some $s,t\in\Zn$. 
Let $\La(c)=s+t=k$ be the level of $\La$. 

\subsection{Paths}
For convenience, we extend the subscript $i$ of $\La_i$ to
$i \in \Z$, by setting $\La_i=\La_{i'}$ for $i\equiv i'\,(\bmod~2)$.
We also set $\widehat{i} = \La_{i+1} - \La_i\,(i = 0,1)$.
For $k\in\Zn$, we define a set of level $k$ weights $P_k$ by
$P_k=\{a_0\La_0+a_1\La_1\mid a_0,a_1\in\Z,a_0+a_1=k\}$. 

\begin{definition}[path on $P_k$]
Fixing $k,L\in\Zn$ we define a path $p$ of length $L$ as 
a sequence $p = (p_0,\cdots,p_L,p_{L+1})$, with all $p_i\in P_k$ 
and 
$p_{i+1}-p_i \in \{ k\widehat{0}, (k-1)\widehat{0}+\widehat{1}, \cdots, 
k\widehat{1} \}$. 
\end{definition}

\begin{definition}[set of paths $\Pset$]
For $L\in\Zn$ and a dominant weight $\La=s\La_0+t\La_1$ $(s+t=k)$, we set
\[
\Pset=\left\{p=(p_0,\cdots,p_L,p_{L+1})\,\left\vert\,
\begin{array}{l}
p \mbox{ is a path on } P_k,\\
p_i=s\La_i+t\La_{i+1}\mbox{ for }i=L,L+1
\end{array}
\right.\right\}.
\]
\end{definition}
To the set $\Pset$, we associate a special path $\bar p$, called 
the ground-state path, defined as follows: 

\begin{definition}[ground-state path $\bar{p}$ ]\label{gsp}
\[
\bar{p}=(\bar{p}_0,\cdots,\bar{p}_L,\bar{p}_{L+1}),\qquad 
\bar{p}_i=s\La_i+t\La_{i+1}\mbox{ for all }i.
\]
\end{definition}
\noindent The relevance of the ground-state path will become
clear in the sequel.

For $k\in\Zn$ set 
$\A^+_k=\{a_0\ep_0+a_1\ep_1\mid a_0,a_1\in\Zn,a_0+a_1=k\}$.
Here we consider $\ep_0$ and $\ep_1$ as symbols. Just as in the
case of $\La_i$, we extend the subscript $i$ of $\ep_i$ to $i\in\Z$ 
be setting $\ep_i = \ep_{i'}$ for $i \equiv i'\,(\bmod2)$.
We encode a path in terms of a sequence of elements in $\A^+_k$ 
as follows:

\begin{definition}[sequence of elements $\iota(p)$]\label{defiota}
Let $\La=s\La_0+t\La_1\,(s+t=k)$.
For a path $p=(p_0,\cdots,p_L,$ $p_{L+1})\in \Pset$ we define 
a sequence $\iota(p)=(\iota(p)_0,\cdots,\iota(p)_L)$ by 
$\iota(p)_i=m_i\ep_0+(k-m_i)\ep_1\in\A^+_k$, where for each 
$i$, $m_i$ is determined by 
$p_{i+1}-p_i=m_i\widehat{0}+(k-m_i)\widehat{1}$.
\end{definition}
For the ground-state path $\bar{p}$, $\iota(\bar{p})$ is given by 
$\iota(\bar{p})_i=s\ep_i+t\ep_{i+1}$. We remark that the data $\iota(p)$
uniquely determine $p\in\Pset$.

Define a function $\mu:\A^+_k\longrightarrow\{0,1,\cdots,k\}$ by 

$$
\mu(m\ep_0+(k-m)\ep_1)=m.
$$
Using $\mu$, we define the energy, and the weight of a path $p$ 
as follows:

\begin{definition}[energy of a path $E(p)$]\label{energy}
\begin{eqnarray*}
E(p)&=&
\sum_{j=1}^L
j\,\big(H(\iota(p)_{j-1},\iota(p)_j)-
           H(\iota(\bar{p})_{j-1},\iota(\bar{p})_j)\big),
\end{eqnarray*}
where 
\begin{equation}\label{defH} 
H(\ep,\ep')=\max(k-\mu(\ep),\mu(\ep'))\quad\mbox{for }\ep,\ep'\in\A^+_k.
\end{equation} 
\end{definition}

\begin{definition}[weight of a path $\wt p$]
\begin{eqnarray*}
\wt p&=&p_0-E(p)\delta.
\end{eqnarray*}
\end{definition}
We note that for the ground-state path $\bar{p}\in\Pset$ we have
$E(\bar{p})=0$ and $\wt \bar{p}=\La$. 
We will need the following lemma in the sequel.
\begin{lemma}\label{gse}
Let $\La=s\La_0+t\La_1\,(s+t=k)$. For the ground-state path
$\bar{p}\in\Pset$, we have
\[
\sum_{j=1}^L jH(\iota(\bar{p}_{j-1}),\iota(\bar{p}_j))
=\left(\frac{L+\ep^{(L)}}2\right)^2 k+(-1)^{\ep^{(L)}}
\frac{L+\ep^{(L)}}2 s.
\]
Here $\ep^{(L)}$ is defined by
\begin{equation} 
\ep^{(L)}=\left\{
\begin{array}{ll}
0&(L:\,\mbox{even}),\\
1&(L:\,\mbox{odd}).
\end{array}\right.\label{defep}
\end{equation} 
\end{lemma}

\Proof
Note that $\iota(\bar{p}_j)=s\ep_j+t\ep_{j+1}$. From (\ref{defH})
we obtain the following expression to calculate:
\[
\sum_{j=1\atop j:\mbox{\tiny odd}}^L jt
+\sum_{j=1\atop j:\mbox{\tiny even}}^L js,
\]
which is evaluated easily.
\qed

So far, we fixed the length of paths $L$. Next, we consider two lengths
$L$ and $L'$ ($L\le L'$). Let $\La=s\La_0+t\La_1$. Then there is an injection 
given by
\begin{eqnarray*}
\Pset&\longrightarrow&\P_{L'}(\La)\\
p=(p_0,\cdots,p_L,p_{L+1})
&\mapsto&
p'=(p'_0,\cdots,p'_{L'},p'_{L'+1}),
\end{eqnarray*}
where $p'_i=p_i$ ($0\le i\le L+1$), $=s\La_i+t\La_{i+1}$ ($L+1\le i\le L'+1$).
Let $\bar{p}'$ be the ground-state path of $\P_{L'}(\La)$. The image of
$\Pset$ in $\P_{L'}(\La)$ is characterized by the paths $\{p'\}$ such that
$p'_i=\bar{p}'$ ($i\ge L$). Note that the energy and weight of a path are
unchanged by this injection. We have the following inductive system:

\[
\P_0(\La)\longrightarrow\P_1(\La)\longrightarrow\cdots
\longrightarrow\Pset\longrightarrow\cdots.
\]
We define the set of all paths $\P(\La)$ as follows:

\begin{definition}[$\P(\La)$]
\[
\P(\La)=\lim_\rightarrow \Pset.
\]
\end{definition}

\subsection{Extended Young diagrams}

An extended Young diagram $Y$ is an infinite sequence $(y_j)_{j\ge0}$
of integers, such that $y_j\le y_{j+1}$ for all $j$ and $y_j=y_\infty$
(a fixed integer) for sufficiently large $j$ \cite{DJKMO,MM,JMMO}. In 
other words, the sequence stabilizes after a finite (though arbitrarily 
large) number of elements. The number $y_\infty$ is called the `charge' 
of $Y$. In this paper, we only consider extended Young diagrams of 
charge $\gamma=0$ or $1$. 

We can view an extended Young diagram $Y=(y_j)_{j\ge0}$ as a diagram 
drawn on the lattice in the right-half plane with sites 
$\{(m,n)\in \Z\times\Z\mid m\ge0\}$ where $y_j$ denotes the 
depth of the $j$-th column. If $y_j\ne y_{j+1}$, for some $j$, then 
there will be concave (`$\lceil$'), and convex (`$\rfloor$')
corners. A corner located at site $(m,n)$ has diagonal number 
$d = m + n$. 

\subsubsection{A coloring scheme}

We assign each corners of $Y=(y_j)_{j\ge0}$ one of two colors: 
$0=$ `white', or $1=$ `black' as follows: A $d$-diagonal corner 
is white (resp. black), if $d$ is even (resp. odd). A white 
(resp. black) corner will also be called a 0-corner (resp. 1-corner). 

Thus, we can also view an extended Young diagram $Y=(y_j)_{j\ge0}$ of 
charge 0 (resp. 1) as a usual Young diagram, but with its nodes 
alternately colored `white' and `black', such that the top left-most 
node has color `white' (resp. `black'). 

We define the weight $\wt Y$ of an extended Young diagram $Y$ of charge 
$i$ ($=0$ or $1$) to be 

\[
\wt Y=\La_i-k_0\alpha_0-k_1\alpha_1,
\]
where $k_0$ is the number of white nodes in $Y$, and $k_1$ is the
number of of black nodes. We also define the `width of $Y$',
$|Y|$, to be the number of nodes in the first row of $Y$. 
Equivalently, $|Y|=L$ if and only if $y_L=y_\infty$(= charge of $Y$) 
but $y_{L-1}<y_L$. 

\subsubsection{Example}

Consider the extended Young diagram 
$Y=(y_j)_{j\ge0}=(-2,-1,-1,0,0,1,1,1,\cdots)$ of charge 1. 
$Y$ has convex 1-corners at sites $(1,-2)$ and $(5,0)$, convex 0-corner 
at site $(3,-1)$, concave 1-corner at site $(3,0)$ and a concave 0-corner 
at site $(1,-1)$. Also note that $\wt Y=\La_1-4\alpha_0-5\alpha_1$ and 
$|Y|=5$.

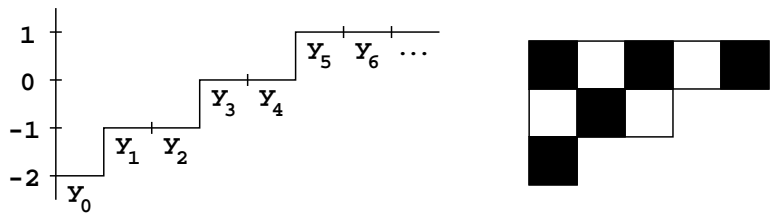
\begin{figure}[h]
\begin{minipage}{7cm}
\centerline{\epsfbox{fig1.eps}}
\end{minipage}
\begin{minipage}{4cm}
\centerline{\epsfbox{fig2.eps}}
\end{minipage}
\end{figure}

\subsection{From paths to extended Young diagrams}

Let us define a `pattern' to be a map
\begin{eqnarray*}
t:\Z\times\Zn & \longrightarrow & \Z\\
(i,j) & \longmapsto & t_{ij}
\end{eqnarray*}
such that
\begin{itemize}
\item[(i)] for all $i$, $(t_{ij})_{j\ge0}$ is an extended Young diagram,
\item[(ii)] $t_{ij}\le t_{i+1\,j}$ for all $i$ and $j$,
\item[(iii)] $t_{i+k\,j}=t_{ij}+2$ for all $i$ and $j$.
\end{itemize}
We say the pattern $t$ is normalized if $0\le\gamma_1\le\cdots\le
\gamma_k<2$, where $\gamma_i=t_{i\infty}$ is the charge of $(t_{ij})_{j\ge0}$.
We call $\gamma=(\gamma_1,\cdots,\gamma_k)$ the charge of $t$. 
Viewing the condition (iii) we can identify the pattern $t$ with a $k$-tuple
$\Y=(Y_1,\cdots,Y_k)$ of extended Young diagrams $Y_i=(t_{ij})_{j\ge0}$.
In terms of extended Young diagrams, (ii) is equivalent to the following
inclusion rule:
\begin{equation} \label{incl} 
Y_1\supset Y_2\supset \cdots\supset Y_k\supset Y_1[2]
\end{equation} 
where $Y_1[2]$ denotes the extended Young diagram which is obtained by
shifting $Y_1$ upward two units on the lattice in the right-half plane.
Let $\La=s\La_0+t\La_1$ ($s+t=k$). Denote the set of patterns with 
charge $\gamma=(\underbrace{0,\cdots,0}_s,\underbrace{1,\cdots,1}_t)$
by $\T(\La)$. Then we have a map
\begin{eqnarray*}
\pi:\T(\La) & \longrightarrow & \P(\La)\\
t=(t_{ij}) & \longmapsto & p,
\end{eqnarray*}
where $p$ is determined from 
$\iota(p)_j=\ep_{t_{1j}+j}+\cdots+\ep_{t_{kj}+j}$.
For a path $p\in\P(\La)$ we say $t$ is a `lift' of $p$ if $t\in\pi^{-1}(p)$.
Among lifts of $p\in\P(\La)$ there exists a unique normalized lift 
$t=(t_{ij})\in\T(\La)$ such that $t_{ij}\ge t'_{ij}$ for all $i,j$ for any
$t'=(t'_{ij})\in\pi^{-1}(p)$ (Proposition 3.4 in \cite{JMMO}).
$t(p)$ is called the `highest lift' of $p$. Let $\Y=(Y_1,\cdots,Y_k)$ be 
the $k$-tuple of extended Young diagrams corresponding to the highest lift 
of $p$. Then we have (Theorem 5.7 in \cite{DJKMO})
\[
\wt p=\wt Y_1+\cdots+\wt Y_k.
\]

\subsection{The crystal $B(\La)$}

Let $(L(\la),B(\La))$  be the crystal base of $V(\la)$. In this subsection
we give the combinatorial rules used to construct the 
crystal $B(\La)$, following the work of \cite{JMMO,MM}.

\subsubsection{The rules}

Let $\La=s\La_0+t\La_1\in P^+$ of level $k=s+t$. We order the elements of 
$\Z\times \{1,2,\cdots,k\}$ as follows. For $(d,j),(d',j')\in
\Z\times \{1,2,\cdots,k\}$ we say 
\[
(d,j)>(d',j')\mbox{ if and only if }d>d'\mbox{ or }d=d'\mbox{ and }j<j'.
\]
Let $\Yset$ denote the set of $k$-tuples of extended Young diagrams 
$\Y=(Y_1,\cdots,Y_s,$ $Y_{s+1},\cdots,Y_k)$ such that $Y_j$ has charge 0
(resp. 1) for $1\le j\le s$ (resp. $s+1\le j\le k$). For $\Y\in\Yset$
and fixed $i$ ($=0$ or $1$), we define the $i$-signature of $\Y$ to be the
sequence $\vep=(\vep_1,\vep_2,\cdots,\vep_m)$\footnote{Note that
we distinguish $\vep_j$ from $\ep_j$ which has appeared in Definition
\ref{energy}.} such that 
\begin{itemize}
\setlength{\itemsep}{0pt}
\item[A1] $\sum_{j=1}^k\sharp\{i$-corners of $Y_j\}=m$,
\item[A2] each $\vep_r$ is either 0 or 1,
\item[A3] we can define $j(r)$ ($1\le j(r)\le k$) and $d(r)$ in such a way
          that $Y_{j(r)}$ has a $d(r)$-diagonal $i$-corner,
\item[A4] if $\vep_r=0$ (resp. $1$) then $Y_{j(r)}$ has a $d(r)$-diagonal
          concave (resp. convex) $i$-corner,
\item[A5] if $r_1<r_2$, then $(d(r_1),j(r_1))>(d(r_2),j(r_2))$.
\end{itemize}

For fixed $i$-signature $\vep=(\vep_1,\vep_2,\cdots,\vep_m)$ we partition
the set $\{1,2,\cdots,m\}=J\sqcup K_1\sqcup \cdots\sqcup K_r$ into 
disjoint subsets by the following procedure:
\begin{itemize}
\setlength{\itemsep}{0pt}
\item[B1] if there is no $j$ such that $(\vep_j,\vep_{j+1})=(0,1)$ define
          $J=\{1,2,\cdots,m\}$,
\item[B2] if there is some $j$ such that $(\vep_j,\vep_{j+1})=(0,1)$ define
          $K_1=\{j,j+1\}$,
\item[B3] apply B1 and B2 above to $\{1,2,\cdots,m\}\setminus K_1$ to choose
          $J$ or $K_2$ and repeat this as necessary to choose $J$ and 
          $K_1,\cdots,K_r$.
\end{itemize}
Let $\vep_J=(\vep_{j_1},\cdots,\vep_{j_r})$, where $J=\{j_1,\cdots,j_r\}$
and $j_1<\cdots<j_r$. We call a 0 or 1 in the $i$-signature `relevant' if
and only if it is in $\vep_J$.

For fixed $i\in\{0,1\}$ and $\Y,\Y'\in\Yset$ suppose that the following 
conditions hold:
\begin{itemize}
\setlength{\itemsep}{0pt}
\item[C1] $\vep=(\vep_1,\cdots,\vep_m)$ and $\vep'=(\vep'_1,\cdots,\vep'_m)$
          are the $i$-signatures of $\Y$ and $\Y'$ respectively,
\item[C2] the partition $\{1,2,\cdots,m\}=J\sqcup K_1\sqcup \cdots\sqcup K_r$
          is same for both $\vep$ and $\vep'$,
\item[C3] there exists $l\in J$ such that $\vep_l=0,\vep'_l=1$, and 
          $\vep_j=\vep'_j=1$ (resp. $0$) if $j\in J$ and $j<l$ (resp. $j>l$).
\end{itemize}
We define $\ft{i}\Y=\Y'$ and $\et{i}\Y'=\Y$ if and only if the above 
conditions hold. If there exists no such $\Y'$ (resp. $\Y$) for $\Y$
(resp. $\Y'$), we define $\ft{i}\Y=0$ (resp. $\et{i}\Y'=0$).

\subsubsection{Example}

Let $\La=\La_0+\La_1$ and take $\Y=(Y_1,Y_2)\in\Yset$, where
$Y_1=(-2,-2,-1,-1,0,0,\cdots)$ is of charge 0 and 
$Y_2=(-1,0,0,1,1,\cdots)$ is of charge 1. \par

\pagebreak
\begin{figure}[htb]
\[
\mbox{\bf Y}=(Y_1,Y_2)=\left(
\begin{minipage}{1.28cm}
\centerline{\epsfbox{fig3.eps}}
\end{minipage},
\begin{minipage}{0.96cm}
\centerline{\epsfbox{fig4.eps}}
\end{minipage}
\right).
\]
\end{figure}

\noindent
The 0-signature of $\Y$ is $\vep=(0,0,1,1,0)$. So using rule B $\vep_J=(0)$.
Hence by using rule C we have \par

\begin{figure}[htb]
\[
\ft0\Y=\Y'=\left(
\begin{minipage}{1.28cm}
\centerline{\epsfbox{fig5.eps}}
\end{minipage},
\begin{minipage}{0.96cm}
\centerline{\epsfbox{fig4.eps}}
\end{minipage}
\right)
\quad\mbox{and}\quad \et0\Y=0.
\]
\end{figure}

\noindent
However, the 1-signature of $\Y$ is $\vep=(1,1,0,0,0)$ and by rule B 
$\vep_J=\vep$. So by rule C we have  \par

\begin{figure}[htb]
\[
\ft1\Y= \left(
\begin{minipage}{1.28cm}
\centerline{\epsfbox{fig6.eps}}
\end{minipage},
\begin{minipage}{0.96cm}
\centerline{\epsfbox{fig4.eps}}
\end{minipage}
\right)
\quad\mbox{and}\quad \et1\Y=\left(
\begin{minipage}{1.28cm}
\centerline{\epsfbox{fig3.eps}}
\end{minipage},
\begin{minipage}{0.96cm}
\centerline{\epsfbox{fig7.eps}}
\end{minipage}
\right).
\]
\end{figure}

\noindent
Thus, in general $\et{i},\ft{i}(i=0,1)$ are well-defined maps from 
$\Yset$ to $\Yset\sqcup\{0\}$. Now let $\Phi=(\phi_1,\cdots,\phi_s,
\phi_{s+1},\cdots,\phi_k)\in \Yset$ be the $k$-tuple of empty extended Young
diagrams where $\phi_j$ has charge 0 (resp. 1) for $1\le j\le s$
(resp. $s+1\le j\le k$). Note that $\wt \Phi=\La=s\La_0+t\La_1$,
and $\et{i}\Phi=0$ for $i=0,1$. We call $\Phi$ the highest weight vector 
(or vacuum) of weight $\La$.

We now define the set $B(\La)\subset\Yset$ as:
\[
B(\La)=\{\ft{i_1}\ft{i_2}\cdots\ft{i_k}\Phi\mid k\ge0,i_j=0,1\}\setminus\{0\}.
\]
For $\Y\in B(\La),i=0,1$, define 
\begin{eqnarray*}
\vep_i(\Y)&=&\max\{m\ge0\mid \et{i}^m\Y\in B(\La)\},\\
\varphi_i(\Y)&=&\max\{m\ge0\mid \ft{i}^m\Y\in B(\La)\}.
\end{eqnarray*}
Thus $\vep_i$\footnote{We use this new definition of $\vep_i$ only here.
We hope this will not cause any confusion.},$\varphi_i:B(\La)
\longrightarrow\Z$. We also have the map
$\mbox{\sl wt}:B(\La)\longrightarrow P$ given by $\Y\mapsto\wt \Y$. 
The set $B(\La)$
equipped with these maps $\et{i},\ft{i},\vep_i,\varphi_i$ and $\wt$ 
turns out to be the `crystal' (see \cite{K2,KKM}) associated with 
the integrable highest weight module $V(\La)$. In particular, for $i=0,1$
we have 
\[
\et{i}(B(\La))\subseteq B(\La)\sqcup\{0\}\quad \mbox{and}\quad
\ft{i}(B(\La))\subseteq B(\La)\sqcup\{0\}.
\]
Moreover, $B(\La)$ turns out to be the set of $k$-tuples of extended 
Young diagrams corresponding to the highest lifts of $\P(\La)$
(Proposition 3.12 in \cite{JMMO}).

The crystal $B(\La)$ has the structure 
of an oriented colored graph with elements of $B(\La)$ as the set of 
vertices and for $\Y,\Y'\in B(\La),\Y\stackrel{i}{\longrightarrow}\Y'
(i=0,1)$ if and only if $\ft{i}\Y=\Y'$. As a graph $B(\La)$ is connected.
In particular, for $\Y\in B(\La)$ we have (see Proposition 3.11 in 
\cite{JMMO}),
\begin{equation} 
\Y\neq\Phi \Longrightarrow \et0\Y\neq0\quad\mbox{or}\quad\et1\Y\neq0.
\label{notann}
\end{equation} 
As a consequence of (\ref{incl}), we have
\[
|Y_1|\ge\cdots\ge|Y_s|\quad\mbox{and}\quad|Y_{s+1}|\ge\cdots\ge|Y_k|.
\]
Now we have the following simple but important observations.

\begin{proposition}\label{noteq}
Let $\Y=(Y_1,\cdots,Y_s,Y_{s+1},\cdots,Y_k)\in B(\La)$ and $\Y\neq\Phi$.
Then $|Y_1|\neq|Y_{s+1}|$.
\end{proposition}

\Proof
Suppose $|Y_1|=|Y_{s+1}|$. Since $Y_1$ has charge 0 and $Y_{s+1}$ has
charge 1, the last nodes in the first rows of $Y_1$ and $Y_{s+1}$ will
have opposite color. Also by (\ref{incl}), the second rows of $Y_1$ and
$Y_{s+1}$ have at least one node less than their first rows. Consequently,
we have two convex corners of opposite colors corresponding to the last
nodes in the first rows of $Y_1$ and $Y_{s+1}$. One of these convex 
corners will contribute a 1 to the 0-signature of $\Y$ and the other
will contribute a 1 to the 1-signature of $\Y$. However, in both cases
this 1 will not be relevent as it will be preceeded by at least one 0
in each signature. Therefore, by applying $\et0$ and $\et1$ successively
we will obtain $\Y'=(Y'_1,\cdots,Y'_s,Y'_{s+1},\cdots,Y'_k)\in B(\La)$ 
such that
\begin{eqnarray*}
&&\Y'=\et{i_m}\et{i_{m-1}}\cdots\et{i_1}\Y
\mbox{ for some }m\mbox{ and }i_j\in\{0,1\},\\
&&\et0\Y'=0,\et1\Y'=0,\\
&&|Y'_1|=|Y_1|=|Y_{s+1}|=|Y'_{s+1}|.
\end{eqnarray*}
Since $\Y\neq\Phi$ and $|Y'_1|=|Y_1|$, we also have $\Y'\neq\Phi$.
So by (\ref{notann}) either $\et0\Y'\neq0$ or $\et1\Y'\neq0$, which
is a contradiction. Hence $|Y_1|\neq|Y_{s+1}|$.
\qed

\begin{proposition}\label{bounded}
Let $\Y=(Y_1,\cdots,Y_s,Y_{s+1},\cdots,Y_k)\in\Yset$. Suppose that
$\ft{i}^n\Y=\Y'=(Y'_1,\cdots,Y'_s,Y'_{s+1},\cdots,Y'_k)$ for some $n>0$
and $i=0$ or $1$. Then $|Y'_j|\le |Y_j|+1$ for all $1\le j\le k$.
\end{proposition}

\Proof
Let $\vep=(\vep_1,\vep_2,\cdots,\vep_m)$ be the $i$-signature of $\Y$
and $\vep_J=(\vep_{j_1},\cdots,\vep_{j_r})$ be the relevant part of 
$\vep$. By definition of $\ft{i}$ action (see rule C), $\ft{i}^n\Y=
\Y'\neq0$ implies that $\vep_{j_{r-n+1}}=\vep_{j_{r-n+2}}=\cdots
=\vep_{j_r}=0$. Also recall that each relevant 0 in $\vep_J$ corresponds
to a unique concave $i$-corner in some $Y_j,1\le j\le k$ and each 
application of $\ft{i}$ changes only one relevant 0 to a relevant 1 
without affecting the partition of $\vep$. Hence the first row of each 
$Y'_j$ can have at most one more $i$-color node than that of $Y_j$ for
$1\le j\le k$. Therefore, the proposition follows.
\qed

\section{Crystals of Demazure modules}

\subsection{Demazure modules}

The affine Lie algebra $\slth$ has a triangular decomposition 
$\en^-\oplus \ha\oplus \en^+$, where $\ha=\mbox{ span }\{h_0,h_1,d\}$ 
is the Cartan subalgebra and $\en^+$ (resp. $\en^-$) denote the sum of
the positive (resp. negative) root spaces (see \cite{Kac}). The subalgebra 
$\bi=\ha+\en^+$ is the Borel subalgebra of $\slth$. Let $r_0,r_1$ denote
the simple reflections corresponding to the simple roots $\alpha_0,\alpha_1$
respectively. Recall that $r_i\mu=\mu-\mu(h_i)\alpha_i$ for all $\mu\in P$.
Let $W$ denote the Weyl group of $\slth$ generated by $r_0,r_1$. For 
$w\in W$ let $l(w)$ denote the length of $w$ and let $\prec$ denote the 
Bruhat order on $W$. For $\La=s\La_0+t\La_1\in P^+,k=s+t$, as before we 
consider the integrable highest weight $\Uq{\slth}$-module $V(\La)$. 
It is known that
for $w\in W$, the extremal weight space $V(\La)_{w\La}$ is one dimensional.
Let $E_w(\La)$ denote the $\Uq{\bi}$-module generated by $V(\La)_{w\La}$.
These modules $E_w(\La),w\in W$ are called the Demazure modules. They
are finite-dimensional subspaces of $V(\La)$ and have the following 
property:
\begin{eqnarray*}
&&\mbox{For }w,w'\in W,w\preceq w',\mbox{ we have }
E_w(\La)\subseteq E_{w'}(\La),\\
&&\mbox{and }\bigcup_{w\in W}E_w(\La)=V(\La).
\end{eqnarray*}

\subsection{Demazure crystals}

It is known that for $\slth$, the Weyl group is
\[
W=\{(r_1r_0)^m,r_0(r_1r_0)^m,(r_0r_1)^m,r_1(r_0r_1)^m\mid m\in\Zn\}.
\]
In particular, for each integer $L>0$, the Weyl group $W$ has two distinct
elements of length $L$:
\[
w^+_L=\underbrace{\cdots r_0r_1r_0}_L\quad\mbox{and}\quad
w^-_L=\underbrace{\cdots r_1r_0r_1}_L.
\]
Let $W^+=\{1,w^+_L\mid L>0\}$ and $W^-=\{1,w^-_L\mid L>0\}$. Note that on 
each set $W^+$ and $W^-$ the Bruhat order is a total order. Also note that
$W^+\cup W^-=W$ and $W^+\cap W^-=\{1\}$. For $\La=s\La_0+t\La_1\in P^+$,
Kashiwara \cite{K2} has defined the crystals $B_w(\La)$ for the Demazure
modules $E_w(\La)$ as a suitable subset of the crystal $B(\La)$. For our 
purpose it suffices to recall the following recursive property of $B_w(\La)$:
\begin{eqnarray} 
&&\mbox{If }r_i w\prec w, \mbox{then}\nonumber\\
&&B_w(\La)=\{\ft{i}^m b\mid m\ge0,b\in B_{r_iw}(\La),\et{i}b=0\}
\setminus\{0\},
\label{rec}
\end{eqnarray} 
which can be used to construct $B_w(\La)$ starting from the vacuum 
$\Phi\in B(\La)$.

\begin{figure}[bp]
\begin{minipage}{11cm}
\centerline{\epsfbox{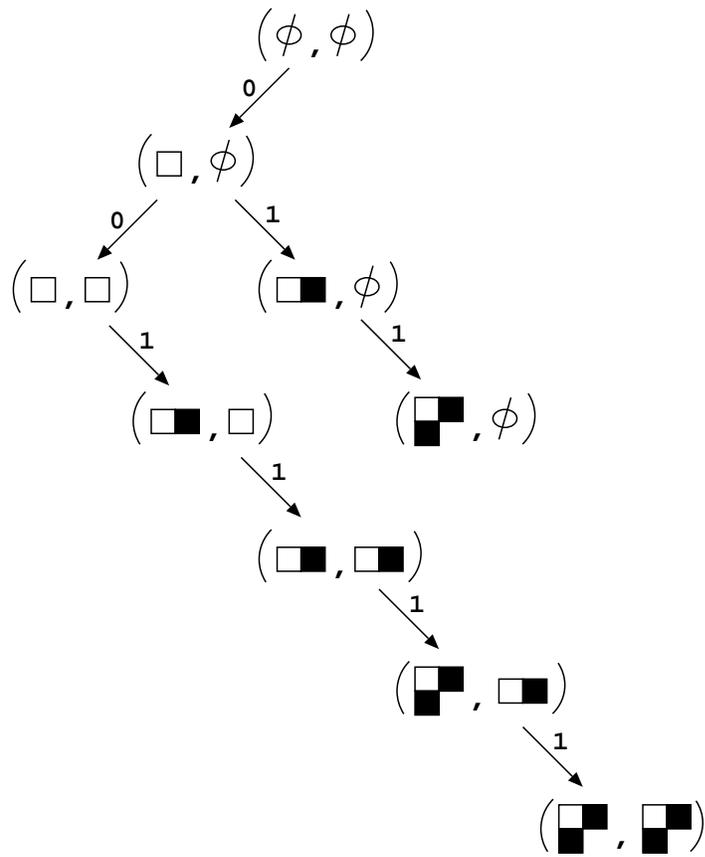}}
\end{minipage}
\caption{The crystal graph of $B_{r_1r_0}(2\La_0)$.}
\label{fig8}
\end{figure}

For example, if $\La=2\La_0$ and $w=r_1r_0$, the graph of the crystal 
$B_w(\La)$ can be easily seen in Figure \ref{fig8}.

Since $w^+_{L-1}\prec w^+_L,w^-_{L-1}\prec w^-_L,
w^+_{L-1}\prec w^+_{L-1} r_1=w^-_L$ and $w^-_{L-1}\prec w^-_{L-1} r_0=w^+_L$,
the following result is an immediate consequence of Proposition 3.2.4 in
\cite{K2}.

\begin{proposition}\label{Demincl}
For $L>0$ and $\La\in P^+$ we have
\begin{itemize}
\setlength{\itemsep}{0pt}
\item[(i)] $B_{w^+_{L-1}}(\La)\subseteq B_{w^+_L}(\La)$,
\item[(ii)] $B_{w^-_{L-1}}(\La)\subseteq B_{w^-_L}(\La)$,
\item[(iii)] $B_{w^+_{L-1}}(\La)\subseteq B_{w^-_L}(\La)$,
\item[(iv)] $B_{w^-_{L-1}}(\La)\subseteq B_{w^+_L}(\La)$.
\end{itemize}
\end{proposition}

\subsection{Realizations of $B_{w^+_L}(\La)$ and $B_{w^-_L}(\La)$}

We define an extended Young diagram $Y=(y_j)_{j\ge0}$ of charge $\gamma$
($=$0 or 1) to be {\em maximal of width} $L$ if 
\[
|Y|=L\mbox{ and }y_{j+1}=y_j+1\mbox{ for }0\le j\le L-1.
\]
We also define, for $\La=s\La_0+t\La_1\in P^+,k=s+t$,
\[
B_L(\La)=\{\Y=(Y_1,\cdots,Y_s,Y_{s+1},\cdots,Y_k)\in B(\La)
\mid |Y_j|\le L,1\le j\le k\}.
\]
Noting that the weight multiplicities of $w^+_L\La$ and $w^-_L\La$ are
one, we define extremal vectors as follows:

\begin{definition}[extremal vectors]
$b_{w^+_L\La}$ (resp. $b_{w^-_L\La}$) is defined as the element of
$B(\La)$ which has the weight $w^+_L\La$ (resp. $w^-_L\La$).
\end{definition}
The following theorem characterizes the extremal vectors.

\begin{theorem}
Let $\La=s\La_0+t\La_1,k=s+t$ and $L>0$. Then we have
\begin{itemize}
\item[(i)] $b_{w^+_L\La}=(Y_1,\cdots,Y_s,Y_{s+1},\cdots,Y_k)\in B(\La)$
where $Y_1=\cdots=Y_s$ are maximal of width $L$ and charge 0, and 
$Y_{s+1}=\cdots=Y_k$ are maximal of width $L-1$ and charge 1,
\item[(ii)] $b_{w^-_L\La}=(Y_1,\cdots,Y_s,Y_{s+1},\cdots,Y_k)\in B(\La)$
where $Y_1=\cdots=Y_s$ are maximal of width $L-1$ and charge 0, and 
$Y_{s+1}=\cdots=Y_k$ are maximal of width $L$ and charge 1.
\end{itemize}
\end{theorem}

\Proof
We use induction on $L$. If $L=1$, then $w^+_L=r_0$ and $w^+_L\La
=r_0(s\La_0+t\La_1)=\La-s\alpha_0$. Also it is easy to see that the 
0-signature of the vacuum $\Phi\in B(\La)$ is $\vep=(\vep_1,\cdots,\vep_s)$
where $\vep_j=0$ for $1\le j\le s$. Hence
\[
b_{r_0\La}=\ft0^s\Phi=(\Box,\cdots,\Box,\phi,\cdots,\phi)\in B(\La)
\]
and (i) holds for $L=1$. Now assume (i) for $L-1$. Observe that 
$b_{w^+_{L-1}\La}$ has convex and concave corners of opposite colors.
Furthermore, $b_{w^+_{L-1}\La}$ has exactly $sL+t(L-1)$ number of concave 
$i$-corners where $L-1\equiv i\,(\bmod2)$. Also
\begin{eqnarray*}
w^+_L\La&=&r_i(w^+_{L-1}\La)=w^+_{L-1}\La-(w^+_{L-1}\La)(h_i)\alpha_i\\
&=&w^+_{L-1}\La-(sL+t(L-1))\alpha_i.
\end{eqnarray*}
Since applying $\ft{i}$ to $b_{w^+_{L-1}\La}$ $(sL+t(L-1))$ times we add one 
$i$-color box to each row of each diagram in $b_{w^+_{L-1}\La}$, we have
\[
b_{w^+_L\La}=\ft{i}^{sL+t(L-1)}b_{w^+_{L-1}\La}
=(Y_1,\cdots,Y_s,Y_{s+1},\cdots,Y_k)\in B(\La)
\]
where $Y_1=\cdots=Y_s$ are maximal of width $L$ and $Y_{s+1}=\cdots=Y_k$
are maximal of width $L-1$. This proves (i). Proof of (ii) is similar
to that of (i).
\qed

The following theorem characterizes the crystals of the Demazure modules.

\begin{theorem}\label{Demcry}
Let $\La=s\La_0+t\La_1,k=s+t$ and $L>0$. Then we have 
\begin{eqnarray*}
B_{w^+_L}(\La)=&&\hspace{-1.8em}\{(Y_1,\cdots,Y_s,Y_{s+1},\cdots,Y_k)\in B(\La)
\mid |Y_1|\le L, |Y_{s+1}|\le L-1\},\\
B_{w^-_L}(\La)=&&\hspace{-1.8em}\{(Y_1,\cdots,Y_s,Y_{s+1},\cdots,Y_k)\in B(\La)
\mid |Y_1|\le L-1, |Y_{s+1}|\le L\}.
\end{eqnarray*}
\end{theorem}

\Proof
Let 
\begin{eqnarray*}
S_L&=&\{(Y_1,\cdots,Y_s,Y_{s+1},\cdots,Y_k)\in B(\La)
\mid |Y_1|\le L, |Y_{s+1}|\le L-1\},\\
T_L&=&\{(Y_1,\cdots,Y_s,Y_{s+1},\cdots,Y_k)\in B(\La)
\mid |Y_1|\le L-1, |Y_{s+1}|\le L\}.
\end{eqnarray*}
We use induction on $L$. First we will show that $B_{w^+_L}(\La)
\subseteq S_L$ and $B_{w^-_L}(\La)\subseteq T_L$. This is clear for 
$L=1$ since by (\ref{rec}) $B_{r_0}(\La)=\{\ft0^m\Phi\mid 0\le m\le s\}$
and $B_{r_1}(\La)=\{\ft1^m\Phi\mid 0\le m\le t\}$. Assume that 
$B_{w^+_{L-1}}(\La)\subseteq S_{L-1}$ and 
$B_{w^-_{L-1}}(\La)\subseteq T_{L-1}$.
Since $w^+_L=r_i w^+_{L-1}$ where $L-1\equiv i\,(\bmod2)$ by (\ref{rec}) 
we have
\[
B_{w^+_L}(\La)=\{\ft{i}^m b\mid m\ge0,b\in B_{w^+_{L-1}}(\La),\et{i}b=0\}
\setminus\{0\}.
\]
Since $B_{w^+_{L-1}}(\La)\subseteq S_{L-1}$, by Proposition \ref{bounded}
we have $B_{w^+_L}(\La)\subseteq S_L$. Similarly $B_{w^-_L}(\La)\subseteq T_L$.

Now we proceed to show that $S_L\subseteq B_{w^+_L}(\La)$ and
$T_L\subseteq B_{w^-_L}(\La)$. As before this is clear for $L=1$. Assume that
$S_{L-1}\subseteq B_{w^+_{L-1}}(\La)$ and $T_{L-1}\subseteq B_{w^-_{L-1}}(\La)$. Let $\Y=(Y_1,\cdots,Y_s,Y_{s+1},\cdots,Y_k)\in S_L$. By Proposition 
\ref{noteq}, $|Y_1|\neq|Y_{s+1}|$. If $|Y_1|\le L-1$ and $|Y_{s+1}|\le L-2$,
then $\Y\in S_{L-1}\subseteq B_{w^+_{L-1}}(\La)\subseteq B_{w^+_L}(\La)$
by Proposition \ref{Demincl}(i).
If $|Y_1|\le L-2$ and $|Y_{s+1}|\le L-1$, then $\Y\in T_{L-1}
\subseteq B_{w^-_{L-1}}(\La)\subseteq B_{w^+_L}(\La)$ by Proposition
\ref{Demincl}(iv). 

Finally, if $|Y_1|=L$ and $|Y_{s+1}|\le L-1$, then $Y_1$ will have a convex
$i$-corner corresponding to the last node of the first row of $Y_1$ where
$L-1\equiv i\,(\bmod2)$. Furthermore, this convex $i$-corner will contribute
a relevant 1 to the $i$-signature of $\Y$. Therefore, $\et{i}\Y\neq0$.
Now, we can choose $m$ such that $\et{i}^m\Y=\Y'=(Y'_1,\cdots,Y'_s,
Y'_{s+1},\cdots,Y'_k)\in B(\La)$ and $\et{i}\Y'=0$. Then $|Y'_1|=L-1$,
and $|Y'_{s+1}|\le L-2$ by Proposition \ref{noteq}. Thus
$\Y'\in S_{L-1}\subseteq B_{w^+_{L-1}}(\La)$. Hence by (\ref{rec})
$\Y=\ft{i}^m\Y'\in B_{r_iw^+_{L-1}}(\La)=B_{w^+_L}(\La)$. This proves that
$S_L\subseteq B_{w^+_L}(\La)$. Similarly $T_L\subseteq B_{w^-_L}(\La)$
and the theorem holds.
\qed

The following result is an immediate consequence of Theorem \ref{Demcry}.

\begin{corollary}
For $\La=s\La_0+t\La_1$, we have 
\begin{eqnarray*}
\mbox{(i)}& &B_{w^+_L}(\La)\cup B_{w^-_L}(\La)=B_L(\La),\\
\mbox{(ii)}& &B_{w^+_L}(\La)\cap B_{w^-_L}(\La)=B_{L-1}(\La).
\end{eqnarray*}
\end{corollary}

Let us consider the extreme cases $s=0$ or $t=0$. Assume $\La=k\La_0$.
Noting $w^-_L\La=w^+_{L-1}\La$ we have $B_{w^-_L}(\La)=B_{w^+_{L-1}}(\La)$.
Similarly, if $\La=k\La_1$, we have $B_{w^+_L}(\La)=B_{w^-_{L-1}}(\La)$.
Combining Proposition \ref{Demincl} and the previous corollary, we have

\begin{corollary} \label{Demextreme}
\[
B_L(\La)=\left\{
\begin{array}{ll}
B_{w^+_L}(\La) & \quad\mbox{if }\La=k\La_0,\\
B_{w^-_L}(\La) & \quad\mbox{if }\La=k\La_1.
\end{array}
\right.
\]
\end{corollary}

\section{Characters}

\subsection{1D configuration sum}
We summarize the 1D configuration sum which have already appeared in the 
study of exactly solvable SOS models \cite{DJKMO}. It will be used in the 
evaluation of the character for $\Pset$ in the following subsection.

Fixing a positive integer $k$, let us call a pair of intergers $(b,c)$ 
{\em weakly admissible} if
$b-c=-k,-k+2,\cdots,k$. For a weakly admissible pair $(b,c)$ and $L\in\Zn$,
let $f^{(k)}_L(b,c)$ denote the unique solution to the following linear 
difference equation.
\begin{eqnarray} 
f^{(k)}_L(b,c)&=&\mathop{{\sum}'}_d f^{(k)}_{L-1}(d,b)q^{L|d-c|/4},
\label{diffeq}\\
f^{(k)}_0(b,c)&=&\delta_{b0}. \label{initial}
\end{eqnarray} 
Here the sum $\sum'$ is taken over $d$ such that the pair $(d,b)$ is 
weakly admissible. If $(b,c)$ is not weakly admissible, $f^{(k)}_L(b,c)$ is
set to be 0. Note that $N,m$ in \cite{DJKMO} are replaced with $k,L$, 
respectively.
$f^{(k)}_L(b,c)$ enjoys the following properties.
$$
\begin{array}{lrcl}
\mbox{Reflection symmetry:}&&&\\
&f^{(k)}_L(b,c)&=&f^{(k)}_L(-b,-c),\\
\mbox{Support property:}&&&\\
&f^{(k)}_L(b,c)&=&0 \quad\mbox{ unless }|b|\le Lk,|c|\le (L+1)k \\
&&&\qquad\qquad\quad\mbox{ and }b\equiv Lk \bmod2.
\end{array}
$$

Given an integer $\mu$ set
\[
R_\mu=\{b\in\Z\mid (\mu-1)k\le b\le (\mu+1)k\}, 
\]
where the left (resp. right) equality sign is taken if $\mu-1\le0$ 
(resp. $\mu+1\ge0$). We also set $R_{\mu,\nu}=R_\mu\times R_\nu$. For
a weakly admissible pair $(b,c)$ and integer $L(\ge1)$, let 
$\mu,\nu(=\mu\pm1)$ be integers such that 
\[
(b,c)\in R_{\mu,\nu},\quad \mu\equiv L+1,\nu\equiv L+2 \bmod2.
\]
The above equation uniquely determines $\mu$ and $\nu$ except when 
$b=0(\mu=\pm1)$ or $c=0(\nu=\pm1)$, in which cases either choice is allowed.

Let us recall standard notations in $q$ analysis (see \cite{Andrews} 
for example). We set
\[
(z;q)_m=\left\{
\begin{array}{cl}
{\displaystyle \prod_{j=1}^m(1-zq^{j-1})}&(m\ge1),\\
1&(m=0).
\end{array}
\right.
\]
If there is no danger of confusion, we abbreviate it as $(z)_m$. The
{\em $q$-multinomial coefficient} is defined as
\begin{equation} \label{multcoeff} 
{M \brack m_1\cdots m_n}_q=\left\{
\begin{array}{cl}
{\displaystyle \frac{(q)_M}{\prod_{j=1}^n(q)_{m_j}}}
&{\mbox{if }m_1,\cdots,m_n\ge0 \atop \qquad\mbox{ and }m_1+\cdots+m_n=M},\\
0&\qquad\;\;\mbox{otherwise}.
\end{array}
\right.
\end{equation} 
${M \brack m_1\,m_2}_q$ is often written as ${M \brack m_1}_q$ 
and called the {\em Gaussian polynomial}. In both cases, if the subscript
is omitted, we understand it as $q$. 

For later use we extend the definition of the Gaussian polynomial
${M \brack i}$ with $M\ge0,i\in\Z$ to the case of $M,i\in\Z$, by
setting
\begin{equation} \label{newdefgau} 
{M \brack i}=\left\{
\begin{array}{cl}
{\displaystyle \frac{(q^{M-i+1})_i}{(q)_i}}
&\qquad\;\;\mbox{if }i\ge0,\\
0&\qquad\;\;\mbox{otherwise}.
\end{array}
\right.
\end{equation} 
Similarly, for $M<0$ we define $(z)_M$ to be $(q^Mz)_{-M}^{-1}$.

\begin{lemma}
For $M,N,n\in\Z$ we have 
\begin{eqnarray} 
(z)_M&=&\sum_i(-z)^i q^{i(i-1)/2}{M \brack i}, \label{gau1}\\
{M+N \brack n}&=&\sum_i q^{(n-i)(M-i)}{M \brack i}{N \brack n-i}, \label{gau2}
\end{eqnarray} 
where the sums are taken over all integers.
\end{lemma}

\Proof
The first equation reduces to (3.3.6) or (3.3.7) in \cite{Andrews} depending
on whether $M\ge0$ or $M<0$.

Next, note that the following equation holds for $M,N\in\Z$:
\begin{equation} 
(z)_{M+N}=(z)_M(q^Mz)_N. \label{forgau2}
\end{equation} 
The second equation can be proved by using (\ref{gau1}) in the both
sides of (\ref{forgau2}) and comparing the coefficients of $z^n$.
\qed

We call the following expression bosonic.

\begin{proposition}[Theorem 4.1.1 in \cite{DJKMO}] \label{DJKMOf}
For any weakly admissible pair $(b,c)\in R_{\mu,\nu}$ and integers $L,k\ge1$,
\begin{eqnarray} 
\lefteqn{(q)_{L-1}f^{(k)}_L(b,c)}\nonumber\\
&=&\left(
\sum_{i\ge(L+\nu)/2 \atop j\le(L+\mu-1)/2}
-\sum_{i\le(L+\nu)/2-1 \atop j\ge(L+\mu+1)/2}\right)
(-1)^{i+j}q^{Q^{(L)}_{i,j}(b,c)}
{L-1 \brack i}{L \brack j}, \label{bos1}\\
\lefteqn{Q^{(L)}_{i,j}(b,c)}\nonumber\\
&=&\frac12(i-j)(i-j+1)-\left(i-\frac{L-1}2\right)\left(j-\frac{L}2\right)k
\nonumber\\
&&+\frac{b}2\left(i-\frac{L-1}2\right)+\frac{c}2\left(j-\frac{L}2\right).\nonumber
\end{eqnarray} 
\end{proposition}

\subsection{Fermionic expression}

In this subsection we rewrite the bosonic expression to the fermionic
expression which involves $q$-multinomials. At the first step we show 
the following.

\begin{proposition} \label{intprop1}
For any weakly admissible pair $(b,c)\in R_{\mu,\nu}$ satisfying $b\ge -Lk$
and integers $L,k\ge1$,
\begin{eqnarray*}
&&\hspace{-3em}q^{\frac{L(L-1)k+(L-1)b+Lc}4}f^{(k)}_L(b,c)\\
&=&\sum_{j\in \Zs}(-1)^jq^{\frac12j(j-1)+\frac{L-1}2jk+\frac{c}2j}
{L \brack j}{L-j-1-(j-\frac{L}2)k+\frac{b}2 \brack L-1}\\
&&-\sum_{j\in \Zs}(-1)^jq^{\frac12j(j+1)+\frac{L}2jk+\frac{b}2j}
{L-1 \brack j}{L-j-2-(j-\frac{L-1}2)k+\frac{c}2 \brack L-1}\\
&&\hspace{15em}\times\left(1-q^{L-j-1-(j-\frac{L-1}2)k+\frac{c}2}\right).
\end{eqnarray*}
\end{proposition}

\Proof
Rewrite the sum
\[
\left(\sum_{i\ge(L+\nu)/2 \atop j\le(L+\mu-1)/2}
-\sum_{i\le(L+\nu)/2-1 \atop j\ge(L+\mu+1)/2}\right)
\]
in (\ref{bos1}) as
\[
\left(\sum_{i\in\Zss \atop j\le(L+\mu-1)/2}
-\sum_{i\le(L+\nu)/2-1 \atop j\in\Zss}\right).
\]
Applying (\ref{gau1}) for the sum over $i\in\Z$ or $j\in\Z$, we obtain
\begin{eqnarray*}
&&\hspace{-3em}q^{\frac{L(L-1)k+(L-1)b+Lc}4}f^{(k)}_L(b,c)\\
&=&\sum_{j\le\frac{L+\mu-1}2}(-1)^jq^{\frac12j(j-1)+\frac{L-1}2jk+\frac{c}2j}
{L \brack j}\frac{(q^{-j+1-(j-\frac{L}2)k+\frac{b}2})_{L-1}}{(q)_{L-1}}\\
&&-\sum_{i\le\frac{L+\nu}2-1}(-1)^iq^{\frac12i(i+1)+\frac{L}2ik+\frac{b}2i}
{L-1 \brack i}\frac{(q^{-i-(i-\frac{L-1}2)k+\frac{c}2})_L}{(q)_{L-1}}.
\end{eqnarray*}
Recall that we have assumed $b\ge-Lk$. Noticing $(\mu-1)k\le b\le (\mu+1)k,
(\nu-1)k\le c\le (\nu+1)k$ with $\mu\ge1-L,\nu\ge-L$, we get 
$L-j-1-(j-\frac{L}2)k+\frac{b}2<m-1$ and 
$L-i-2-(i-\frac{L-1}2)k+\frac{c}2<m-1$ if $j>\frac{L+\mu-1}2$ and 
$i>\frac{L+\nu}2-1$. Therefore, recalling the definition (\ref{newdefgau}), 
we have
\begin{eqnarray*}
{L-j-1-(j-\frac{L}2)k+\frac{b}2 \brack L-1}
&=&
\left\{
\begin{array}{ll}
\frac{(q^{-j+1-(j-\frac{L}2)k+\frac{b}2})_{L-1}}{(q)_{L-1}}
&\quad \mbox{if }j\le\frac{L+\mu-1}2\\
0&\quad \mbox{otherwise},
\end{array}
\right.\\
{L-i-2-(i-\frac{L-1}2)k+\frac{c}2 \brack L-1}
&=&
\left\{
\begin{array}{ll}
\frac{(q^{-i-(i-\frac{L-1}2)k+\frac{c}2})_{L-1}}{(q)_{L-1}}
&\quad \mbox{if }i\le\frac{L+\nu}2-1\\
0&\quad \mbox{otherwise}.
\end{array}
\right.
\end{eqnarray*}
Thus we arrive at the desired result.
\qed

\begin{proposition} \label{prop6}
For any weakly admissible pair $(b,c)\in R_{\mu,\nu}$ satisfying 
$b\ge0,c\neq b+k$ and integers $L\ge0,k\ge2$, we have
\begin{eqnarray*}
f^{(k)}_L(b,c)=\sum_{0\le i\le L}
&&\hspace{-1.5em}q^{\frac{L(L-1)}4-\frac{(k-1)i^2+(2L+b+c-1)i}4}{L \brack i}\\
&&\times f^{(k-1)}_{L-i}(b+(k+1)i-L,c+(k+1)i-L+1).
\end{eqnarray*}
\end{proposition}

\Proof
If $L=0$, the formula is shown easily from the initial condition 
(\ref{initial}). We assume $L\ge1$.

Notice that from the assumption $f^{(k-1)}_{L-i}(b+(k+1)i-L,c+(k+1)i-L+1)=0$
when $i=L$. Since the right hand side of the formula in Proposition 
\ref{intprop1}
is formally zero when $L=0$, we can use it for convenience of our proof.
Check that the assumption allows us to use Proposition \ref{intprop1}. Setting
\begin{eqnarray*}
a^{(k)}_L(b,c;j)&=&
q^{-\frac{L(L-1)k+(L-1)b+Lc}4+\frac{L-1}2jk+\frac{c}2j}
{L \brack j}{L-j-1-(j-\frac{L}2)k+\frac{b}2 \brack L-1},\\
b^{(k)}_L(b,c;j)&=&
q^{-\frac{L(L-1)k+(L-1)b+Lc}4+\frac{L}2jk+\frac{b}2j}
{L-1 \brack j}\\
&&\times
{L-j-2-(j-\frac{L-1}2)k+\frac{c}2 \brack L-1}
\left(1-q^{L-j-1-(j-\frac{L-1}2)k+\frac{c}2}\right),
\end{eqnarray*}
we have to show
\begin{eqnarray*}
&&\hspace{-3em}\sum_{j\in\Zs}(-1)^j
\left(q^{\frac12j(j-1)}a^{(k)}_L(b,c;j)
-q^{\frac12j(j+1)}b^{(k)}_L(b,c;j)\right)\\
&=&\sum_{i=0}^L q^{\frac{L(L-1)}4-\frac{(k-1)i^2+(2L+b+c-1)i}4}
{L \brack i}\\
&&\hspace{1.4em}\times\sum_{j\in\Zs}(-1)^j
\left(q^{\frac12j(j-1)}a^{(k-1)}_{L-i}(b',c';j)
-q^{\frac12j(j+1)}b^{(k-1)}_{L-i}(b',c';j)\right),
\end{eqnarray*}
where $b'=b+(k+1)i-L,c'=c+(k+1)i-L+1$.
In fact, we can show that the recurrence relation holds term by term,
that is,
\begin{eqnarray} 
a^{(k)}_L(b,c;j)&=&\sum_{i=0}^L
q^{\frac{L(L-1)}4-\frac{(k-1)i^2+(2L+b+c-1)i}4}{L \brack i}\nonumber\\
&&\hspace{1.4em}\times a^{(k-1)}_{L-i}(b+(k+1)i-L,c+(k+1)i-L+1;j), \label{recofa}
\end{eqnarray} 
and the same for $b^{(k)}_L(b,c;j)$.
After some calculation, (\ref{recofa}) reduces to 
\begin{eqnarray*}
&&\hspace{-3em}
{L \brack j}{L-j-1-(j-\frac{L}2)k+\frac{b}2 \brack L-1}\\
&=&\sum_i q^{(L-i-1)(L-i-j)}{L \brack i}
{L-i \brack j}{-1-(j-\frac{L}2)k+\frac{b}2 \brack L-i-1}.
\end{eqnarray*}
Noting that ${L \brack i}{L-i \brack j}={L \brack j}{L-j \brack i}$,
we can cancel the factor ${L \brack j}$. Thus we arrive at (\ref{gau2})
with $M=L-j,N=-1-(j-\frac{L}2)k+\frac{b}2,n=L-1$.

Similarly, (\ref{recofa}) with $a^{(k)}_L(b,c;j)$ replaced by $b^{(k)}_L(b,c;j)$
reduces to
\begin{eqnarray} 
&&\hspace{-3em}
{L-1 \brack j}{L-j-2-(j-\frac{L-1}2)k+\frac{c}2 \brack L-1}
\left(1-q^{L-j-1-(j-\frac{L-1}2)k+\frac{c}2}\right)\nonumber\\ 
&=&\sum_i q^{(L-i)(L-i-j-1)}{L \brack i}
{L-i-1 \brack j}{-1-(j-\frac{L-1}2)k+\frac{c}2 \brack L-i-1}\nonumber\\
&&\hspace{15em}\times\left(1-q^{-(j-\frac{L-1}2)k+\frac{c}2}\right).\label{recofb}
\end{eqnarray} 
Set $N=-(j-\frac{L-1}2)k+\frac{c}2$. Noting that
\begin{eqnarray*}
{L-1 \brack j}{L-j+N-2 \brack L-1}\left(1-q^{L-j+N-1}\right)\\
&&\hspace{-10em}={L \brack j}{L-j+N-1 \brack L}\left(1-q^{L-j}\right),\\
{L \brack i}{L-i-1 \brack j}{N-1 \brack L-i-1}\left(1-q^N\right)\\
&&\hspace{-10em}={L \brack j}{N \brack L-i}{L-j-1 \brack i}\left(1-q^{L-j}\right),
\end{eqnarray*}
and cancelling the factor ${L \brack j}\left(1-q^{L-j}\right)$, (\ref{recofb})
reduces again to (\ref{gau2}) with $M=L-j-1,n=L$.
\qed

This recurrence relation with respect to $k$ leads us to an expression
of $f^{(k)}_L(b,c)$ in terms of the $q$-multinomial coefficient 
(\ref{multcoeff}).
Next theorem is a generalization of Theorem 4.2.5 in \cite{DJKMO}.

\begin{theorem} \label{th3}
For any weakly admissible pair $(b,c)\in R_{\mu,\nu}$ and integers
$L\ge0,k\ge1$, we have the following expression for $f^{(k)}_L(b,c)$.
\begin{eqnarray} 
f^{(k)}_L(b,c)&=&\sum q^{\cal Q}{L \brack x_0 \cdots x_k}, \label{th3eq1}\\
{\cal Q}&=&\frac{L^2k-L(b-c+k)+b}4-\sum_{0\le a\le a'\le k}(a'-a)x_ax_{a'}\nonumber\\
&&\hspace{8em}+\sum_{\frac{c-b+k}2\le a\le k}(a-\frac{c-b+k}2)x_a.\nonumber
\end{eqnarray} 
The sum in (\ref{th3eq1}) is taken over all non negative integers
$x_0,\cdots,x_k$ satisfying
\[
\sum_{a=0}^k x_a=L,\qquad
\sum_{a=0}^k ax_a=\frac{Lk-b}2.
\]
\end{theorem}

\Proof
We prove by induction on $k$. If $k=1$, (\ref{th3eq1}) reduces to 
\[
f^{(1)}_L(b,c)=q^{\frac{bc}4}{L \brack \frac{L+b}2\,\frac{L-b}2}.
\]
This trivially satisfies (\ref{initial}). Checking (\ref{diffeq}) reduces 
to the following identity of Gaussian polynomials.
\[
{M \brack i}={M-1 \brack i-1}+q^i{M-1 \brack i}
=q^{M-i}{M-1 \brack i-1}+{M-1 \brack i}.
\]
Thus the theorem is correct when $k=1$.

Now assume $k\ge2$. If $c=b+k$, (\ref{th3eq1}) turns out to be equivalent to 
Theorem 4.2.5 (see also equation (4.2.1)) in \cite{DJKMO}. Therefore we can 
assume $c\neq b+k$. First suppose $b\ge0$. From Proposition \ref{prop6} and 
the assumption of the induction, we obtain
\begin{equation} \label{th3eq2} 
f^{(k)}_L(b,c)=\sum_{x_k=0}^L\sum q^{\cal Q}
{L \brack x_k}{L-x_k \brack x_0\cdots x_{k-1}},
\end{equation} 
with
\begin{eqnarray*}
{\cal Q}&=&\frac{L(L-1)}4-\frac{(k-1)x_k^2+(2L+b+c-1)x_k}4+{\cal Q'},\\
{\cal Q'}&=&\frac{(L-x_k)^2(k-1)-(L-x_k)(b'-c'+k-1)+b'}4\\
&&-\sum_{0\le a<a'\le k-1}(a'-a)x_ax_{a'}
+\sum_{\frac{c'-b'+k-1}2\le a\le k-1}(a-\frac{c'-b'+k-1}2)x_a,\\
b'&=&b+(k+1)x_k-L,\quad c'=c+(k+1)x_k-L+1,
\end{eqnarray*}
and the inner sum in (\ref{th3eq2}) is taken over non negative integers
$x_0,\cdots,x_{k-1}$ such that $\sum_{a=0}^{k-1}x_a=L-x_k$,
$\sum_{a=0}^{k-1}ax_a=\frac{(L-x_k)(k-1)-b'}2$.
Calculating the power of $q$ and noting that
\[
{L \brack x_k}{L-x_k \brack x_0 \cdots x_{k-1}}
={L \brack x_0 \cdots x_k},
\]
we arrive at (\ref{th3eq1}).

For the proof of the case $b<0$, it suffices to check Reflection symmetry 
in the right hand side of (\ref{th3eq1}). The check is easily done by 
the variable change $x_j\rightarrow x_{k-j}$ ($0\le j\le k$) in the 
summand of (\ref{th3eq1}).
\qed

\subsection{Character of $\Pset$}
We define the character for $\Pset$. 
\begin{definition}[character of $\Pset$] \label{defch}
\[
\ch_L(\La)=\ch_L(\La)(z,q)=\sum_{p\in\Pset}z^{(\La-p_0\mid\La_1)}q^{E(p)}.
\]
\end{definition}

Then we have the following expression for $\ch_L(\La)$.

\begin{proposition}\label{pathch}
\[
\ch_L(\La)=\sum_{j\in\Zs}z^{-j}q^{j/2}
f^{(k)}_L(\ep^{(L)}(s-t)-2j,\ep^{(L+1)}(s-t)-2j).
\]
Here $\ep^{(L)}$ is given in (\ref{defep}).
\end{proposition}

\Proof
First, for $a\La_0+b\La_1\in P_k$, we set $\nu(a\La_0+b\La_1)=b$. Then
we have
\[
\mu(\iota(p)_j)=\frac{\nu(p_{j+1})-\nu(p_j)+k}2.
\]
From the definition of $H$ (Definition \ref{energy}), we have
\begin{eqnarray*}
H(\iota(p)_{j-1},\iota(p)_j)
&=&\frac{k}2+\frac{|\nu(p_{j-1})-\nu(p_{j+1})|}4\nonumber\\
&&+\frac{\nu(p_{j-1})-2\nu(p_j)+\nu(p_{j+1})}4,
\end{eqnarray*}
and therefore,
\begin{eqnarray*}
\sum_{j=1}^L jH(\iota(p)_{j-1},\iota(p)_j)
&=&\frac{k}2\frac{L(L+1)}2
+\sum_{j=1}^L j\frac{|\nu(p_{j-1})-\nu(p_{j+1})|}4\nonumber\\
&&+\frac{\nu(p_0)-(L+1)\nu(p_L)+L\nu(p_{L+1})}4.
\end{eqnarray*}
Thus we have 
\[
E(p)=\frac{\nu(p_0)-\nu(\bar{p}_0)}4
+\sum_{j=1}^L j\frac{|\nu(p_{j-1})-\nu(p_j)|}4.
\]
Using this and parametrizing $p_0$ by $\La+j\alpha_1$ ($j\in\Z$), we can 
rewrite the character (Definition \ref{defch}) as
\begin{eqnarray*}
\ch_L(\La)&=&\sum_{j\in\Zs}z^{-j}
q^{(\nu(p_0)-\nu(\bar{p}_0))/4}K_{L,j},\\
K_{L,j}&=&\sum_{p\in\Pset \atop p_0=\La+j\alpha_1}
q^{\sum_{j=1}^L j|\nu(p_{j-1})-\nu(p_{j+1})|/4}.
\end{eqnarray*}
Noting that $K_{L,j}$ satisfies the same difference equation as 
$f^{(k)}_L(\nu(p_L)-\nu(p_0),\\ \nu(p_{L+1})-\nu(p_0))$, we have
\[
\ch_L(\La)=\sum_{j\in\Zs}z^{-j}
q^{(\nu(p_0)-\nu(\bar{p}_0))/4}
f^{(k)}_L(\nu(p_L)-\nu(p_0),\nu(p_{L+1})-\nu(p_0))
\]
with $p_0=\La+j\alpha_1$. Evaluating $\nu(p_j)$ ($j=0,L,L+1$) explicitly,
we have the desired result.
\qed

Now we have the final result on the character of $\Pset$.

\begin{theorem} \label{fermionic}
Let $\La=s\La_0+t\La_1$ ($s+t=k$), $C$ be the Cartan matrix of $sl(k)$,
that is, $C_{ij}=2\delta_{ij}-\delta_{i\,j+1}-\delta_{i\,j-1}$
($1\le i,j\le k-1$). Set $x^t=(x_1,\cdots,x_{k-1})$ and let $e_i$ be
the $k-1$ dimensional $i$-th unit vector ($e_0=e_k=0$). Define
$F_{j,L}(\La)$ by
\[
\ch_L(\La)=\sum_{j\in\Zs}z^{-j}F_{j,L}(\La).
\]
Then we have the following expression for $F_{j,L}(\La)$.
\par\noindent
If $L$ is even,
\[
F_{j,L}(\La)=\sum_{\sum_{i=1}^{k-1}ix_i\in j+k\Zs}
q^{x^tC^{-1}x-x^tC^{-1}e_s+\frac{j(j+t)}k}
{L \brack x_0\cdots x_k},
\]
where
\[
x_0=\frac{L}2-\frac1k\left(\sum_{i=1}^{k-1}(k-i)x_i+j\right),
\quad
x_k=\frac{L}2-\frac1k\left(\sum_{i=1}^{k-1}ix_i-j\right).
\]
If $L$ is odd,
\[
F_{j,L}(\La)=\sum_{\sum_{i=1}^{k-1}ix_i\in j+t+k\Zs}
q^{x^tC^{-1}x-x^tC^{-1}e_t+\frac{j(j+t)}k}
{L \brack x_0\cdots x_k},
\]
where
\[
x_0=\frac{L+1}2-\frac1k\left(\sum_{i=1}^{k-1}(k-i)x_i+j+t\right),
\quad
x_k=\frac{L-1}2-\frac1k\left(\sum_{i=1}^{k-1}ix_i-j-t\right).
\]
\end{theorem}

\Proof
Suppose $L$ is even, from Theorem \ref{th3} and Proposition \ref{pathch},
we obtain
\[
F_{j,L}(\La)=\sum_{\sum_{i=0}^k ix_i=\frac{Lk}2+j}
q^{\frac{L^2k}4-\frac{Lt}2-\sum_{0\le i\le i'\le k}(i'-i)x_ix_{i'}
+\sum_{s\le i\le k}(i-s)x_i}
{L \brack x_0\cdots x_k}.
\]
Solving $\sum_{i=0}^k x_i=L$ and $\sum_{i=0}^k ix_i=\frac{Lk}2+j$ for
$x_0$ and $x_k$, we get 
$x_0=\frac{L}2-\frac1k\left(\sum_{i=1}^{k-1}(k-i)x_i+j\right)$,
$x_k=\frac{L}2-\frac1k\left(\sum_{i=1}^{k-1}ix_i-j\right)$.
Substitute these expressions into the power of $q$, and note that
$(C^{-1})_{ij}=\frac{i(k-j)}k$ (for $i\le j$), $\frac{j(k-i)}k$ (for $i>j$).
Noticing that the sum restriction changes into 
$\sum_{i=1}^{k-1}ix_i\in j+k\Zs$, we obtain the desired expression.

The case of $L$ odd is similar.
\qed

\subsection{Character of Demazure module}
As we have seen in section 2, crystals of Demazure modules can be realized
as subsets of the crystal $B(\La)$ or equivalently $\P(\La)$. 
Now it is easy to see that the characters
of Demazure modules admit the following definition.

\begin{definition}
\[
\ch_L^\pm(\La)=\ch_L^\pm(\La)(z,q)=\sum_{p\in B_{w_L^\pm}(\La)}
z^{(\La-p_0\mid \La_1)}q^{E(p)}.
\]
\end{definition}
A natural question to ask is: Is it possible to express $\ch_L^\pm(\La)$
in terms of $\ch_L(\La)$~? The answer is given by the following theorem.

\begin{theorem}\label{Demch}
Let $\La=s\La_0+t\La_1,s+t=k$. We have 
\begin{eqnarray*}
\ch_L^+(\La)&=&\sum_{i=0}^s z^{-\ep^{(L)}(s-i)}
q^{\frac{L+\ep^{(L)}}2(s-i)} \ch_{L-1}(\La^{(i)}),\\
\ch_L^-(\La)&=&\sum_{i=0}^t z^{\ep^{(L)}(t-i)}
q^{\frac{L-\ep^{(L)}}2(t-i)} \ch_{L-1}(\La^{(k-i)}).
\end{eqnarray*}
Here $\ep^{(L)}$ is given in (\ref{defep}) and 
$\La^{(i)}=i\La_0+(k-i)\La_1$.
\end{theorem}

\Proof
First we prove the expression for $\ch_L^+(\La)$. Let us suppose $L$ is
even. Interpreting Theorem \ref{Demcry} in terms of paths, for any 
$p\in B_{w_L^+}(\La)$ we have $\iota(p)_{L-1}
=(k-i)\ep_0+i\ep_1$ where $i=0,1,\cdots,s$. For each $i$
consider the following subset of $\P_L(\La)$.

\[
S^{(i)}=\{p\in\P_L(\La)\mid p_{L-1}=\La-((k-i)\widehat{0}+i\widehat{1})\}.
\]
Noting that $\La-((k-i)\widehat{0}+i\widehat{1})-\La^{(k-i)}=\La-\La^{(i)}$,
we see that 
\[
\{(p_j-(\La-\La^{(i)}))_{j=0,\cdots,L}\mid p\in S^{(i)}\}
\]
can be identified with $\P_{L-1}(\La^{(i)})$. Defining $\bar{p}^{(i)}$
as the ground-state path of $\P_{L-1}(\La^{(i)})$, we have 
\begin{eqnarray*}
\ch_{L-1}(\La^{(i)})&=&
\sum_{p\in B_{w_L^+}(\La) \atop p_{L-1}=\La-((k-i)\widehat{0}+i\widehat{1})}
z^{(\La^{(i)}-(p_0-(\La-\La^{(i)}))\mid \La_1)}\nonumber\\
&&\qquad\qquad \times 
q^{\sum_{j=1}^{L-1}j(H(\iota(p)_{j-1},\iota(p)_j)
-H(\iota(\bar{p}^{(i)})_{j-1},\iota(\bar(p)^{(i)})_j))}.
\end{eqnarray*}
Therefore,
\begin{eqnarray*}
\ch_L^+(\La)&=&\sum_{i=0}^s q^{\sum_{j=1}^{L-1}j
(H(\iota(\bar{p}^{(i)})_{j-1},\iota(\bar{p}^{(i)})_j)
-H(\iota(\bar{p})_{j-1},\iota(\bar{p})_j))}\nonumber\\
&&\qquad \times 
q^{L(H(\iota_{L-1}^{(i)},\iota_L^{(i)})-H(\iota_{L-1},\iota_L))}
\ch_{L-1}(\La^{(i)}).
\end{eqnarray*}
Here $\iota_{L-1}^{(i)}=(k-i)\ep_0+i\ep_1,
\iota_{L-1}=\iota_{L-1}^{(s)},\iota_L^{(i)}=\iota_L=
s\ep_0+t\ep_1$.
Using Lemma \ref{gse} we get the result. The case of $L$ odd is similar.

For $\ch_L^-(\La)$ we utilize the Dynkin diagram symmetry.
\begin{eqnarray*}
\ch_L^-(\La)&=&\ch_L^+(t\La_0+s\La_1)(qz^{-1},q)\\
&=&\sum_{i=0}^t(qz^{-1})^{-\ep^{(L)}(t-i)}
q^{\frac{L+\ep^{(L)}}2(t-i)}\ch_{L-1}(\La^{(i)})(qz^{-1},q).
\end{eqnarray*}
Noting that $\ch_{L-1}(\La^{(i)})(qz^{-1},q)=\ch_{L-1}(\La^{(k-i)})(z,q)$,
we get the result.
\qed

\subsection{Some specializations}
In this subsection, we consider two specializations of Demazure characters, 
and compare with Sanderson's results \cite{S}.

First we consider `real characters'. These are obtained by setting 
$e^{-\alpha_0}=q,e^{-\alpha_1}=q^{-1}$. We have 

\begin{corollary}
Let $\La=s\La_0+t\La_1,s+t=k$.
\[
\ch_L^+(\La)(q^{-1},1)=q^{-\frac{(L-\ep^{(L)})k}2}[s+1][k+1]^{L-1}.
\]
Here $[a]=(1-q^a)/(1-q)$.
\end{corollary}

\Proof
From Theorem \ref{fermionic} we easily get
\[
f^{(k)}_L(b,c;1)=\sum_{\sum_{i=0}^k ix_i=\frac{Lk-b}2}
{L \choose x_0\cdots x_k}.
\]
Here ${L \choose x_0 \cdots x_k}$ denotes the usual multinomial coefficient.
Assume $L$ is even. From Proposition \ref{pathch} we have
\begin{eqnarray*}
\ch_L(\La)(q^{-1},1)
&=&\sum_{j\in\Zs}q^j\sum_{\sum_{i=0}^k ix_i=\frac{Lk}2+j}
{L \choose x_0\cdots x_k}\\
&=&q^{-\frac{Lk}2}\sum_{x_0,\cdots,x_k}q^{\sum_{i=0}^k ix_i}
{L \choose x_0\cdots x_k}\\
&=&q^{-\frac{Lk}2}(1+q+\cdots+q^k)^L\\
&=&q^{-\frac{Lk}2}[k+1]^L.
\end{eqnarray*}
Similarly, for $L$ odd we have 
\[
\ch_L(\La)(q^{-1},1)=q^{-\frac{Lk-s+t}2}[k+1]^L.
\]
Using Theorem \ref{Demch}, we obtain the desired result.
\qed

Next, we consider `principal characters' obtained by setting 
$e^{-\alpha_0}=e^{-\alpha_1}=q$. In considering this case, we restrict
ourselves to the extreme case $t=0$.

\begin{corollary}
Let $\La=k\La_0$.
\[
\ch_L^+(\La)(q,q^2)=\sum_{x_0,\cdots,x_k}
q^{2x^tC^{-1}x+\frac{k}2{\cal S}({\cal S}+1)}
{L \brack x_0\cdots x_k}_{q^2},
\]
where we have set 
\begin{equation} \label{defS} 
{\cal S}=\frac1k\sum_{i=0}^k(k-2i)x_i.
\end{equation} 
\end{corollary}

\Proof
Note from Corollary \ref{Demextreme} that when $\La=k\La_0$ we have
$\ch_L^+(\La)=\ch_L(\La)$. Using Theorem \ref{fermionic} we obtain
\begin{equation} \label{Dempsc} 
\ch_L^+(\La)(q,q^2)=\sum_{j\in\Zs}q^{-j}
\sum_{\sum_{i=1}^{k-1}ix_i\in j+k\Zs}
q^{2x^tC^{-1}x+\frac{2j^2}k}
{L \brack x_0\cdots x_k}_{q^2},
\end{equation} 
where
\begin{equation} \label{x0xk} 
x_0=\left[\frac{L+1}2\right]-\frac1k\left(\sum_{i=1}^{k-1}(k-i)x_i+j\right),
\quad
x_k=\left[\frac{L}2\right]-\frac1k\left(\sum_{i=1}^{k-1}ix_i-j\right)
\end{equation} 
with $[\,]$ being the Gauss symbol: $[x]$ is the largest integer that does
not exceed $x$. Since the condition 
$\sum_{i=1}^{k-1}ix_i\in j+k\Z$ is equivalent to the condition that all $x_k$ are 
integers, we can rewrite the expression in (\ref{Dempsc}) as
\[
\sum_{x_0,\cdots,x_k}
q^{2x^tC^{-1}x+\frac{2j^2}k-j}
{L \brack x_0\cdots x_k}_{q^2}.
\]
Solving (\ref{x0xk}) with respect to $j$, we get
\[
j=\frac{k}2\left(\left[\frac{L+1}2\right]-\left[\frac{L}2\right]
-{\cal S}\right).
\]
Noting that $\left[\frac{L+1}2\right]-\left[\frac{L}2\right]=\ep^{(L)}$
and substituting this into the expression, we arrive at
\[
\sum_{x_0,\cdots,x_k}
q^{2x^tC^{-1}x+\frac{k}2{\cal S}({\cal S}+1-2\ep^{(L)})}
{L \brack x_0\cdots x_k}_{q^2}.
\]
If $L$ is odd, we can make the variable change $x_i\rightarrow x_{k-i}$,
which amounts to ${\cal S}\rightarrow-{\cal S}$. Thus we get the 
desired result.
\qed

In \cite{S}, Sanderson has evaluated the same character involving the
$q$-multinomial coefficient with argument $q$, as opposed to our expression 
with $q^2$.
Accordingly, we get the following polynomial identity.
\begin{eqnarray*}
\lefteqn{\sum_{x_0,\cdots,x_k} 
q^{2x^tC^{-1}x+\frac{k}2{\cal S}({\cal S}+1)}
{L \brack x_0\cdots x_k}_{q^2}}\\
&&=\sum_{0\le i_1\le\cdots\le i_k\le L}
q^{\frac{i_1(i_1+1)}2+\cdots+\frac{i_k(i_k+1)}2}
{L \brack L-i_k\,i_k-i_{k-1}\cdots i_2-i_1\,i_1}_q,
\end{eqnarray*}
where ${\cal S}$ is defined in (\ref{defS}).

\section{Discussion}

We characterized all Demazure crystals associated with $\slth$ in terms
of paths. We then obtained explicit expressions for their full characters.
Our results extend Sanderson's
results, and reduce to hers when specialized to certain limits. Our
derivations are based on ideas and techniques that originated in
computations of certain physical quantities in exactly solvable lattice
models, particularly Baxter's corner transfer matrix method, and the
combinatorial structures related to it.

Since the very same methods and techniques are available, we expect that
our approach can be used for the cases of other affine algebras too.
Let us consider the $\widehat{sl}(n)$ case. In the case we treated in
this work, the affine Weyl group $W$ had two subsets $W^+$ and $W^-$
satisfying $W^+\cup W^-=W$, $W^+\cap W^-=\{1\}$. On each subset the Bruhat
order was a total order. But if $n>2$ the Weyl group does not have
such a property. Thus it is very difficult to treat the Demazure
modules corresponding to all Weyl group elements. But for some particular type of elements
there does exist a correspondence between Demazure crystals and paths, which 
we would like to report in near future.

As for 1-dimensional configuration sums in the higher rank case, few
results are available except the level 1 case \cite{DJKMO1d}. We would
like to push forward these studies. Along this line we note the works
by Kuniba-Nakanishi-Suzuki \cite{KNS} and Georgiev \cite{Georgiev}. 
In \cite{KNS} they presented conjectures of fermionic expressions
for string functions of vacuum modules $V(k\La_0)$ when the affine
algebra is of type $X_r^{(1)}$. Georgiev proved the conjecture in the 
$A_r^{(1)}$ case. For our purpose we have to 'finitize' their
conjecture or result, since we need generating functions of paths 
of finite length. This would be a challenging problem.

Our approach can be also extended to types other than $\widehat{sl}(n)$.
In these cases we do not have descriptions of crystals by extended Young
diagrams or similar notions yet. We might have to work directly on paths.
We have so far discussed the cases of affine algebra modules. There are
definitely the cases of coset modules left. We also would like to understand
representation theoretical meaning of the truncated generating functions.

\section*{Acknowledgement}

We wish to thank S. Dasmahapatra, A. Kuniba and S.~O.~Warnaar
for discussions and interest in this work.
This work was started while K.~M. and M.~O. were visiting 
the University of Melbourne in July, 1995, and continued while O.~F. was 
visiting the University of Utrecht. We thank the Australian 
Research Council, and the Netherlands Organization for 
Scientific Research (NWO) for their financial support of 
these visits.

\end{document}